\documentclass[a4paper,10pt]{article}
\usepackage{amsmath}
\usepackage{amssymb,amsfonts,amsthm}
\usepackage[utf8]{inputenc} %acentos
\usepackage[english]{babel}

\usepackage{anysize}
\usepackage{enumerate}
\usepackage{tikz}
\usepackage{subfig}
\usepackage[all]{xy}
\usepackage{mathrsfs}
\usepackage{verbatim}
\usepackage{cite}
\definecolor{green}{rgb}{0,0.8,0.5}

\usepackage[pdftex]{hyperref}
\hypersetup{backref,colorlinks=true,linkcolor=blue,citecolor=red}
\usepackage[width=0.8\textwidth]{caption}

\renewenvironment{abstract}{\small\quotation\noindent
 {\bfseries \abstractname .}}{\endquotation \par}

\catcode`\@=11

\def\resetthefootnote{\renewcommand{\thefootnote}{\@arabic\c@footnote} }
\def\@principiremex#1{\trivlist
 \item[\hskip \labelsep{\bfseries #1\ \thethm.}]\ignorespaces}
\def\opar@principiremex#1[#2]{\trivlist
 \item[\hskip \labelsep{\bfseries #1\ \thethm\ (#2).}]\ignorespaces}

\newcommand{\newTHEOremrom}[2]{\newenvironment{#1}{\refstepcounter{thm}\@ifnextchar[{\opar@principiremex{#2}}
{\@principiremex{#2}}}{\qedB\endtrivlist}} \catcode`\@=12
\DeclareMathSymbol{\square}{\mathord}{AMSa}{"03}
\newcommand{\qedB}{\nopagebreak\hspace*{\fill}$\square$\par}
\newcommand{\Qed}{\nopagebreak\hspace*{\fill}{\vrule width6pt height6pt depth0pt}\par}

\newTHEOremrom{defi} {Definition}
\newTHEOremrom {rem} {Remark}
\newTHEOremrom {ex} {Example}

\renewcommand{\geq}{\geqslant}
\renewcommand{\epsilon}{\varepsilon}

\newtheorem{thm}{Theorem}[section]

\newtheorem{lema}[thm]{Lemma}

\newcommand{\R}{\mathbb{R}}

\title{\textbf{The monotonicity of the apsidal angle using the theory of potential oscillators}
}

\author{D. Rojas\\[10pt]
{\small \textsl{Departamento de Matem\'atica Aplicada,}}\\
\vspace{-2pt}
{\small \textsl{Universidad de Granada, Granada, Spain}}\\[5pt]
}

\date{}

\begin{document}

\maketitle

\begin{abstract}
In a central force system the angle between two successive passages of a body through pericenters is called the apsidal angle. In this paper we prove that for central forces of the form $f(r)\sim \lambda r^{-(\alpha+1)}$ with $\alpha<2$ the apsidal angle is a monotonous function of the energy, or equivalently of the orbital eccentricity.
\end{abstract}

\section{Introduction}

The angle between two successive passages of a body in a central force system through pericenters is called the apsidal angle $\Theta$. Its behaviour has attracted both physicists and mathematicians attention since 1687 when Newton stated its precession theorem (Book I, Philosophi\ae\ Naturalis Principia Mathematica): for orbits close to the circular ones, a force proportional to $1/r^{\alpha+1}$ leads to $\Theta=\pi/\sqrt{2-\alpha}$. An important and immediate consequence of this result is that experimental measurements of the apsidal angle close to central orbits may give the exponent of the force law. In 1873 Bertrand~\cite{Bertrand} published a note in Comptes Rendus to prove that among all central field of forces in the Euclidean space there are only two exceptional cases (the harmonic oscillator and the Newtionian potential) in which all solutions close to the circular motions are also periodic. The equivalence of Bertrand's theorem in terms of the apsidal angle is that $\Theta$ remains constant and commensurable with $2\pi$. Indeed this only happens when $\alpha=1$ (the inverse-square Newton's law) and when $\alpha=-2$ (the linear oscillator). In force laws with exponent $\alpha\geq 2$ bounded orbits no longer exist (see~\cite{Danby}) so the study of the apsidal angle can be restricted to $\alpha<2$.

Let us consider a central field of forces in the plane
\[
\ddot x + V'(\|x\|) = 0,\ x\in\R^2
\]
where $V:(0,+\infty)\rightarrow\R$ is a smooth function. It is well-known that the previous system is integrable in the Liouville sense with two first integrals in involution: the energy of the system $h=\frac{1}{2}\|\dot x\|^2+V(\|x\|)+\frac{\ell^2}{2\|x\|^2}$ and the angular momentum $\ell=\|x \wedge \dot x\|$. The apsidal angle of an unbounded orbit with energy $h$ and angular momentum $\ell$ can be written as
\[
\Theta=\int_{r_-}^{r_+}\frac{\ell dr}{r^2\sqrt{2(h-V(r))-\frac{\ell^2}{r^2}}},
\]
where $r=\|x\|$ and $r_{\pm}$ are the distances of the apsis from the center of force and correspond to the solutions of $h-V(r)-\frac{\ell^2}{2r^2}=0$. We point out that circular motions correspond to local minima of $V(r)+\frac{\ell^2}{2r^2}$. In the specific case of a power-law potential of the form $V(r)=-\frac{1}{\alpha}r^{-\alpha}$ the previous equality yield to the function
\begin{equation}\label{apsidal}
\Theta_{\alpha}(h,\ell)=\int_{r_-}^{r_+}\frac{\ell dr}{r^2\sqrt{2(h+\frac{1}{\alpha}r^{-\alpha})-\frac{\ell^2}{r^2}}}.
\end{equation}
Here there is only one radius that yields to circular motion, $r^*=\ell^{\frac{2}{2-\alpha}}$. The previous expression is only valid when the energy $h$ and the angular momentum $\ell$ produce bounded orbits on the system $\ddot r + r^{-\alpha-1}-\ell^2 r^{-3}=0$.

At this point it is common to introduce the orbital eccentricity $e=(r_+-r_-)/(r_+ + r_-)$ in order to study the behaviour of $\Theta_{\alpha}(h,\ell)$ in terms of $e$. Indeed both energy $h$ and angular momentum $\ell$ can be expressed as functions of $r_{\pm}$. This allows the apsidal angle to depend only on $e$.
It should be mentioned that the integral in~\eqref{apsidal} cannot be solved in closed form (except for the cases $\alpha=1$ and $\alpha=-2$, for which remains constant). Expressions in terms of elliptic functions are also given when $\alpha=-1,-\frac{2}{3},\frac{1}{2},\frac{2}{3}$ (see~\cite{Whit}). 
Recently, an analytic proof of the monotonicity of the apsidal angle as a function of the orbital eccentricity has been given by Castelli~\cite{Castelli}, showing that it decreases for any $\alpha\in(-2,1)$. In this work we propose a shorter proof of this monotonicity for the values $\alpha\in(\frac{1}{2},1)$ and we also proof that the apsidal angle is monotonous increasing for any $\alpha\in(1,2)$. On account of the duality given by the relation 
\begin{equation}\label{relation}
\Theta_{\hat\alpha}=\frac{2-\alpha}{\alpha}\Theta_{\alpha}
\end{equation}
with $(2-\alpha)(2-\hat\alpha)=4$ (see~\cite{GR}), we additionally obtain the monotonicity of the apsidal angle for any $\alpha\in(-\infty,-\frac{2}{3})$. Consequently, on account of the result in~\cite{Castelli}, the monotonicity of $\Theta_{\alpha}$ for any $\alpha<2$ is proven. The key point on the proof is the interpretation of the apsidal angle as the period of an abstract oscillator. Such interpretation was also fruitful in~\cite{OR} when giving a shorter proof of Bertrand's theorem. In this occasion we shall use the Schaaf's monotonicity criterion for the time function of potential systems (see Theorem~\ref{Schaaf}).
We prove the following.

\begin{thm}\label{main}
For any $\ell\neq 0$ and any $\alpha\in(-\infty,-\frac{2}{3})\cup(1/2,2)$ the apsidal angle $\Theta_{\alpha}(h)$ is a monotonic function of the energy with $\Theta_{\alpha}'(h)<0$ if $\alpha\in(-2,-\frac{2}{3})\cup(1/2,1)$ and $\Theta_{\alpha}'(h)>0$ if $\alpha\in(-\infty,-2)\cup(1,2)$.
\end{thm}

We point out that for any fixed angular momentum $\ell\neq 0$, the angular eccentricity $e$ is an increasing function of the energy. Indeed using the implicit function theorem with the equality $h+\frac{1}{\alpha}r_{\pm}^{-\alpha}-\frac{\ell^2}{2r_{\pm}^2}=0$ one have that
\[
r_{\pm}'(h)=\frac{r_{\pm}(h)^{\alpha+3}}{r_{\pm}^2-\ell^2 r_{\pm}(h)^{\alpha}},
\]
which implies $r_{-}'(h)<0$ and $r_{+}'(h)>0$ due to $r_-(h)<\left(\frac{1}{\ell^2}\right)^{\frac{1}{\alpha-2}}<r_+(h)$. Thus 
\[
e'(h)=\frac{2\bigl(r_-(h)r_+'(h)-r_-'(h)r_+(h)\bigr)}{\bigl(r_+(h)+r_-(h)\bigr)^2}>0.
\]
Consequently the result in Theorem~\ref{main} can be also stated with the apsidal angle $\Theta_{\alpha}$ depending on the orbital eccentricity, as the author in~\cite{Castelli} does.

\section{Proof of Theorem~\ref{main}}

We will proof Theorem~\ref{main} by using the theory of planar potential oscillators. Consider a potential differential system
$\ddot x +V'(x)=0$ where $V$ is an analytic function on some interval $I$ that contains $x=0$. We suppose $V'(0)=0$ and $V''(0)>0$, so the origin is a non-degenerated center. This means that in the phase space $(x,\dot x)$ there is a neighbourhood of the origin such that all orbits passing through this neighbourhood are closed. The largest neighbourhood with this property is called period annulus and its projection on the $x$-axis will be denoted by $\mathcal J$. 
The main tool we shall use in order to prove Theorem~\ref{main} is the following Schaaf's monotonicity criterion for potential centers (see~\cite{Schaaf}).

\begin{thm}[Schaaf's criterion]\label{Schaaf}
Let $\ddot x + V'(x)=0$ be an analytic potential differential system with a non-degenerated center at the origin and consider its period function $T(h)$ parametrized by the energy. Then $T'(h)>0$ for all $h\in(0,h_0)$ in case that
\begin{enumerate}
\item[$(I_1)$]$5V^{(3)}(x)^2-3V''(x)V^{(4)}(x)>0$ for all $x\in\mathcal J$ with $V''(x)>0$, and
\item[$(I_2)$]$V'(x)V^{(3)}(x)<0$ for all $x\in\mathcal J$ with $V''(x)=0$.
\end{enumerate}
On the other hand, $T'(h)<0$ for all $h\in(0,h_0)$ in case that
\begin{enumerate}
\item[$(D)$]$5V^{(3)}(x)^2-3V''(x)V^{(4)}(x)<0$ for all $x\in\mathcal J$ with $V''(x)\geq 0$.
\end{enumerate}
\end{thm}

The crucial idea is that the apsidal angle can be interpreted as two-dimensional oscillator and so Theorem~\ref{main} follows as a corollary of the previous result. For any $\ell\neq 0$, the differential equation $\ddot r +r^{-\alpha-1}-\ell^2 r^{-3}=0$ produce bounded orbits for $\alpha>0$ when $h\in\bigl(\frac{\alpha-2}{2\alpha}\ell^{\frac{2\alpha}{\alpha-2}},0\bigr)$. Moreover, on account of the symmetry with respect to $\ell=0$ one can restrict to $\ell>0$.
We consider the Clairaut's change of variable $\rho=\frac{\ell}{r}$ on the expression~\eqref{apsidal}. Then the apsidal angle for $\ell>0$ rewrites
\begin{equation}\label{apsidal2}
\Theta_{\alpha}(h)=\int_{\rho_-}^{\rho_+}\frac{d\rho}{\sqrt{2(h-W_{\ell,\alpha}(\rho))}}
\end{equation}
where $W_{\ell,\alpha}(\rho)\!:=\frac{1}{2}\rho^2-\frac{\ell^{-\alpha}}{\alpha}\rho^{\alpha}$ and $\rho_{\pm}=\ell/r_{\mp}$. We point out that expression in~\eqref{apsidal2} coincides with the value of the period of the solution of the oscillator $\ddot\rho + W_{\ell,\alpha}'(\rho)=0$ at the energy level $\frac{\dot\rho^2}{2}+W_{\ell,\alpha}(\rho)=h$.
The key point to apply Schaaf's criterion to the potential differential system $\ddot \rho + W_{\ell,\alpha}'(\rho)=0$ is that, as one can easily verify, we can write the ``test functions" as

\begin{align*}
5W_{\ell,\alpha}^{(3)}(\rho)^2-3W_{\ell,\alpha}''(\rho)W_{\ell,\alpha}^{(4)}(\rho)&=\ell^{-2\alpha}\rho^{2\alpha-6}P_{\alpha}(\ell^{\alpha}\rho^{2-\alpha}),\\
W_{\ell,\alpha}'(\rho)W_{\ell,\alpha}^{(3)}(\rho)&=\ell^{-2\alpha}\rho^{2\alpha-4}Q_{\alpha}(\ell^{\alpha}\rho^{2-\alpha}),\\
W_{\ell,\alpha}''(\rho)&=\ell^{-\alpha}\rho^{\alpha-2}R_{\alpha}(\ell^{\alpha}\rho^{2-\alpha}),
\end{align*}
with 
\begin{align*}
P_{\alpha}(z)&=(\alpha-2)(\alpha-1)^2(2\alpha-1)+3(\alpha-3)(\alpha-2)(\alpha-1)z,\\
Q_{\alpha}(z)&=(\alpha-2)(\alpha-1)-(\alpha-2)(\alpha-1)z,\\
R_{\alpha}(z)&=-(\alpha-1)+z.
\end{align*}
Accordingly we have the following result.

\begin{lema}\label{lema:conditions}
The conditions $(I_1)$, $(I_2)$ and $(D)$ of Schaaf's monotonicity criterion applied to the potential system  $\ddot \rho + W_{\ell,\alpha}'(\rho)=0$ are equivalent to
\begin{enumerate}
\item[$(I_1')$]$P_{\alpha}(z)>0$ for any $z\in\varphi(\mathcal J)$ with $R_{\alpha}(z)>0$,
\item[$(I_2')$]$Q_{\alpha}(z)<0$ for any $z\in\varphi(\mathcal J)$ with $R_{\alpha}(z)=0$,
\item[$(D')$]$P_{\alpha}(z)<0$ for any $z\in\varphi(\mathcal J)$ with $R_{\alpha}(z)\geq 0$,
\end{enumerate}
respectively, where $\varphi(\rho)\!:=\ell^{\alpha}\rho^{2-\alpha}$.
\end{lema}

At this point we recall that on account of the identity~\eqref{relation}, it is enough to proof Theorem~\ref{main} for any $\alpha\in(1/2,2)$.
The projection on the $x$-axis of the period annulus of the center of system $\ddot \rho + W_{\ell,\alpha}'(\rho)=0$ is the interval $\mathcal J=(0,f(\ell,\alpha))$ with $f(\ell,\alpha)=(2\ell^{-\alpha}/\alpha)^{\frac{1}{2-\alpha}}$. Hence $\varphi(\mathcal J)=\bigl(0,\frac{2}{\alpha}\bigr)$. Let us define $\mathcal L_{\alpha}\!:=\{z\in\varphi(\mathcal J): R_{\alpha}(z)>0\}$. If $\alpha\in(1/2,1)$ then $\mathcal L_{\alpha}=\varphi(\mathcal J)$. By Lemma~\ref{lema:conditions} and the previous discussion, condition $(D)$ in Theorem~\ref{Schaaf} is equivalent to requiring that the linear polynomial $P_{\alpha}$ is negative on the interval $\mathcal L_{\alpha}=\bigl(0,\frac{2}{\alpha}\bigr)$. We point out that $P_{\alpha}(0)=(\alpha-2)(\alpha-1)^2(2\alpha-1)<0$ and $P_{\alpha}(2/\alpha)=\alpha^{-1}(\alpha-2)^2(\alpha-1)(9+\alpha+2\alpha^2)<0$ for all $\alpha\in(1/2,1)$. Thus the requirement $(D)$ holds if $\alpha\in(1/2,1)$ and so by Theorem~\ref{Schaaf} the first assertion in Theorem~\ref{main} is proved. On the other hand, if $\alpha\in(1,2)$ then $\mathcal L_{\alpha}=\bigl(\alpha-1,\frac{2}{\alpha}\bigr)$. In this case the linear polynomial $P_{\alpha}$ satisfies $P_{\alpha}(\alpha-1)=5(\alpha-2)^2(\alpha-1)^2>0$ and $P_{\alpha}(2/\alpha)>0$ for all $\alpha\in(1,2)$. Thus by Lemma~\ref{lema:conditions} $(I_1)$ is satisfied if $\alpha\in(1,2)$. Finally, since $R_{\alpha}(z)=0$ if and only if $z=\alpha-1$ and $Q_{\alpha}(\alpha-1)=(1-\alpha)(\alpha-2)^2<0$ for all $\alpha\in(1,2)$ then $(I_2)$ is satisfied in this case. Consequently, by Theorem~\ref{Schaaf} the second assertion in Theorem~\ref{main} holds and the result is proved. 

\begin{rem}
The period function of a slighly different potential of the form $V_{p,q}'(x)=x^p-x^q$, $p,q\in\R$ with $p>q$ was studied in~\cite{MRV2017}. We want to emphasize that the similarity of $W_{\ell,\alpha}'(\rho)=\rho-\ell^{\alpha}\rho^{\alpha-1}$ with $V_{1,\alpha-1}'(x)$ suggest that Theorem~\ref{main} can be interpreted as a corollary of~\cite[Theorem A]{MRV2017}.
\end{rem}

\section*{Acknowledgements}
The author want to thank Prof. Rafael Ortega for the fruitful discussions that led to the interpretation of the apsidal angle as the period function of an abstract oscillator. The author is partially supported by the MINECO grant MTM2014-52209-C2-1-P and MEC/FEDER grant MTM2014-52232-P.
Conflict of interest: The author declare that he have no conflict of interest.

%\bibliography{mybibfile}

\end{document}